\documentclass[journal]{IEEEtran}
\usepackage[cmex10]{amsmath}
\usepackage{dsfont}
\usepackage{graphicx}
\usepackage{adjustbox}
\usepackage{url, multicol, multirow, color, enumerate, verbatim, amsfonts}
\usepackage{threeparttable}
\usepackage{algorithmic}
\usepackage[ruled,vlined]{algorithm2e}
\usepackage{rotating}
\usepackage{amsmath,bm}
\usepackage{tikz}
\usepackage{tabularx}
 \usepackage{stfloats}
 \usepackage{pgfplots}
 
\usetikzlibrary{shapes,arrows,snakes}

\newcommand{\lgreen}{\textcolor{black}}
\newcommand{\fred}{\textcolor{black}}

\newcommand{\lred}{\textcolor{black}}
\newcommand{\zb}{\textcolor{black}}
\newcommand\Bb{{\mathcal {B}}}

\newcommand\Gg{{\mathcal {G}}}

\newcommand\Tt{{\mathcal{T}}}

\newcommand {\Rmnum} [1] {\expandafter\@slowromancap\romannumeral#1\@}

\newcommand{\ignore}[1]{}

\begin{document} \title{Flexibility Management in  Economic Dispatch with Dynamic  Automatic Generation Control}

\markboth{Submitted, {June}~2020}
{Shell \MakeLowercase{\textit{et al.}}: Bare Demo of IEEEtran.cls
for Journals}

\author{Lei Fan,~\IEEEmembership{Senior Member,~IEEE;}
		Chaoyue Zhao,~\IEEEmembership{ Member,~IEEE;}
        Guangyuan Zhang,~\IEEEmembership{Member,~IEEE;}
         Qiuhua  Huang,~\IEEEmembership{Member,~IEEE}

\thanks{Lei Fan is with the Department of Engineering Technology at University of Houston, E-mail: lfan8@uh.edu.  Chaoyue Zhao is with the Department of Industrial and Systems Engineering at University of Washington at Seattle, E-mail:cyzhao@uw.edu. Guangyuan Zhang is a senior commercial Analyst at  RWE Renewables, E-mail: zhangguangyuan1985@gmail.com. Qiuhua Huang is senior power system research engineer at Pacific Northwest National Laboratory, E-mail: qiuhua.huang@pnnl.gov.}
}

\maketitle
\vspace{-2.0cm}
\begin{abstract}
As the \lred{installation of electronically interconnected renewable energy resources grows rapidly in  power systems, system frequency maintenance \zb{and \lgreen{control} become} challenging problems \zb{to maintain the system reliability} in bulk power systems. \zb{As two of the most important frequency control actions in the control centers of independent system operators (ISOs) and utilities,} the interaction between Economic Dispatch (ED) and Automatic Generation Control (AGC) attracts more and more attention}. In this paper, we propose a robust optimization based framework to measure the system flexibility by considering the interaction between two hierarchical processes (i.e.,ED and AGC). We propose a cutting plane algorithm with the reformulation technique to obtain \lred{seven different indices} of the system. In addition, we study the impacts of several system \lred{factors (i.e., the budget of operational cost, ramping capability, and transmission line capacity)} and show numerically how these factors can \lred{influence} the system flexibility.
\end{abstract}

\begin{IEEEkeywords}
 Economic Dispatch, Automatic Generation Control (AGC), Flexibility Management, Robust Optimization, Cutting Plane Method
\end{IEEEkeywords}
\section{Nomenclature}
\addcontentsline{toc}{section}{Nomenclature}
\begin{IEEEdescription}[\IEEEusemathlabelsep\IEEEsetlabelwidth{$V_1,V_2$}]
\item [A. Sets]
\item[${\Bb}$] Set of buses.
\item[$\mathcal{{G}}$] Set of generators.
\item[$\mathcal{{G}}^b$] Set of generators at bus $b$.
\item [$\mathcal{L}$] Set of transmission lines.
\vspace{1mm}

\item[B. Parameters]

\item[$CP_n$] Penalty cost for generator $n$ in dynamic AGC constraints.
\item[$\tau$] Budget for \zb{the} total operational cost.
\item[$F_{l}$] Transmission capacity {of} transmission line $l$ {(MW)}.
\item[${\mbox{P}_n^{\mbox{\tiny max}}}$] Maximum generation amount {(MW)} of generator $n$.
\item[${\mbox{P}_{n}^{\mbox{\tiny min}}}$] Minimum generation amount {(MW)} of generator $n$.
\item[$\mbox{REG}^{\mbox{\tiny U}}_{n}$] Maximum regulation up amount of generator $n$.
\item[$\mbox{REG}^{\mbox{\tiny MinU}}$] Minimum regulation up requirement of the system.
\item[$\mbox{REG}^{\mbox{\tiny D}}_{n}$] Maximum regulation down amount of generator $n$.
\item[$\mbox{REG}^{\mbox{\tiny MinD}}$] Minimum regulation down requirement of the system.
\item[$\mbox{SR}_{n}^{\mbox{\tiny max}}$] Maximum spinning reserve amount of generator $n$.
\item[$\mbox{SR}^{\mbox{\tiny min}}$] Minimum spinning reserve requirement of the system.
\item[$\mbox{RUR}_{n}$] Maximum ramping up rate  of generator $n$.
\item[$\mbox{RDR}_{n}$] Maximum ramping down rate  of generator $n$.
\item[$\mbox{SF}_{b,l}$] Shift factors of transmission line $l$ and bus $b$.
\item[$\bar{d}_{b}$] Nominal load amount at bus $b$.
\item[$\Delta\bar{d}_{t}$] Nominal system load change amount at sub time interval $t$.
\item[$\hat{d}_{b}$] The maximum deviation amount from nominal load amount at bus $b$.
\item[$\Delta\hat{d}_{t}$] The maximum load disturbance amount from the nominal value at sub time interval $t$.
\item[$\Delta{\omega}_{t}^{\mbox{\tiny min}}$] Minimum system frequency change at sub time interval $t$.
\item[$\Delta{\omega}_{t}^{\mbox{\tiny max}}$] Maximum system frequency change at sub time interval $t$.

\vspace{1mm}

\item[C. Random Parameters]
\item[$d_{b}$] \zb{Random} load on bus $b$ (MW).
\item[$\Delta{d}_{t}$] \zb{Random} system load disturbance (MW) at sub time interval $t$.

\vspace{1mm}

\item[D. Decision Variables]
\item[${{oc}_{n}}$] {Generation} {cost} function of generator $n$.
\item[$\lambda^{up}_{b}$] The scale of upper deviation of load at bus $b$.
\item[$\lambda^{dn}_{b}$] The scale of lower deviation of load at bus $b$.
\item[$\lambda^{up}_{t}$] The scale of upper deviation of system load within one sub time interval t.
\item[$\lambda^{dn}_{t}$] The scale of lower deviation of system load within one sub time interval t.
\item[$p_{n}$] Generation amount of generator $n$ (MW).
\item[$\Delta{f}_{n,t}^{\mbox{\tiny GV}+}$] Slack variable for dynamic AGC constraints for generator $n$ at sub time interval $t$.
\item[$\Delta{f}_{n,t}^{\mbox{\tiny GV}-}$] Slack variable for dynamic AGC constraints for generator $n$ at sub time interval $t$.
\item[$reg^{\mbox{\tiny U}}_{n}$] Regulation up amount of generator $n$ (MW).
\item[$reg^{\mbox{\tiny D}}_{n}$] Regulation down amount of generator $n$ (MW).
\item[${sr_{n}}$] Spinning reserve of generator $n$ {(MW)}.
\item[$ \Delta{p}^{\mbox{\tiny GV}}_{n,t}$] The governor power generation change of generator $n$ at sub time interval $t$.
\item[$ \Delta{p}^{\mbox{\tiny M}}_{n,t}$] The prime power generation change of generator $n$ at sub time interval $t$.
\item[$\Delta{\omega}_{t}$] System frequency change.

\end{IEEEdescription}

\section{Introduction} \label{sec: Int}
As the 3D (decarbonization, digitization and decentralization) trend becomes the mainstream in the evolution of energy systems, more and more \lred{renewable energy resources} are installed in the bulk power systems.  For example, as projected by U.S. Energy Information Administration (EIA) \cite{EIA2020AEOL}, the total share of electricity generation from the renewable energy will be $38\%$ by 2050. In particular, solar energy and wind energy will contribute $17.5\%$ and $12.54\%$ of the total electricity generation by 2050 respectively. Moreover, New York state plans to reach $100\%$ \lgreen{carbon-free} by 2050 \cite{NYISO2019CP}, and California state sets \lgreen{its} $100\%$  clean electric power goal by 2045 \cite{CAISO2019CP}. The increase of these carbon-free and non-dispatchable energy resources requires the enhancement of the digital management capability of ISOs and \lred{utilities}.

In order to hedge against the variability and uncertainty of renewable energy generation and therefore achieve a high penetration of renewable energy to the power system, the concept of flexibility has been proposed and investigated, to gauge the capability of the power system in addressing the \lred{variability} of the \lred{net demand (demand net of wind and solar)} \cite{hobbs1992flx,ela2014evolution}. From the time scale's perspective, the flexibilities on planning and operations of the power system have been recently investigated. For example, the index of insufficient ramping resource expectation (IRRE) is proposed in \cite{lannoye2012evaluation} to reflect the flexibility of the power system in the generation expansion \lgreen{planning}. Operational flexibility and local flexibility for the power system operated by transmission system operators are discussed in \cite{bucher2015quantification}. From market design's perspective, ISOs developed different commodity products for the electricity market to capture the flexibility of the resources. For example, CAISO and MISO developed \lgreen{market-based} flexible ramping products that can improve the availability of  the system's ramping capacity \cite{wang2016enhancing}. ERCOT designed a fast frequency response product to maintain sufficient primary frequency control capability \lgreen{under a high} penetration of renewable energy \cite{liu2018participation}. In addition, MISO is investigating short-term reserve products to enhance the system flexibility and ensure reserve deliverability \cite{Wang2019STRP}. All these new \lred{market designs} effectively provide the pricing signal to   flexible resources in the energy market system.  On the resource level, researchers have investigated flexible generation resources such as combined-cycle power plants \cite{fan2015edge} and \cite{dai2018configuration}, pump-storage plants \cite{garcia2008stochastic}, battery energy storage \cite{ding2019rectangle} and \cite{xu2017factoring}. These complex operation models of multi-cycle or multi-stage energy resources in the electricity market not only can strengthen the capability of the grid to respond to the dynamic change of \lred{net demand}, but also can reduce the operational cost of the electricity grid \cite{guan2018unified}. In addition, the modeling approaches for aggregations of flexible resources such as virtual power plants and distributed energy resources aggregators have also been studied in \cite{babaei2019data} and \cite{chen2019aggregate}.

In this paper, we focus on the real-time flexibility of the power system, motivated by the flexibility metric framework proposed by the researchers from ISO-NE \cite{zhao2015unified}. In \cite{zhao2015unified}, the system flexibility is measured in four dimensions (i.e., time, action, uncertainty, and cost). In the time dimension, the short-term flexibility indicates the capability of the system in responding to emergencies or contingencies from minutes to hours. The long-term flexibility indicates the capability of the system in adapting the change of the generation portfolio, system topology, and regulator policy. In the action dimension, different  \lred{control schemes} (e.g., automatic generation control, economic dispatch, unit commitment, outage management, generation and transmission expansion) can be taken by system operators based on \lgreen{different response time windows} to address \lgreen{the variability} of the net demand. In the uncertainty dimension, system operators need to manage  resources to tackle \lgreen{randomnesses} such as equipment failures or forecasting errors of the net demand. Then, the cost restricts the availability of \lred{control schemes}. With these, the overall system \lred{flexibility} of the real-time economic dispatch under the uncertainty can be calculated in a systematic way.

In the current electricity market practice, ED, \lred{which usually runs every 5 minutes,}  provides generation resources with the base dispatch point and \lred{regulation reserve capacity}. However, \lred{the traditional ED process does not consider the impact of dynamic AGC, which \lgreen{is executed} every 2-6 seconds, with the objective of maintaining the system frequency. This conventional ED-AGC hierarchical model may not be able to provide sufficient system flexibility under \zb{a} high renewable penetration to track the second-to-second net demand variation  \cite{li2015miso} and \cite{Zhang2019AGCED}.} Therefore, we will study the real-time flexibility management by explicitly considering the interaction between ED and the dynamic AGC.  
The contributions of our paper can be summarized as follows.
\begin{enumerate}
	\item We propose a robust optimization based framework to measure the system flexibility of the dynamic AGC constrained economic dispatch process. Then, we develop a separation framework to effectively solve the robust feasibility problem. 
	\item We \lred{propose seven flexibility indices for \zb{a} systematic evaluation \zb{of the system flexibility,} and} analyze the system characteristics by using the real-time flexibility management tool to help the system operator understand the impacts of multiple system factors (e.g., budget, ramping capability, and transmission line capability) on the system flexibility.  

\end{enumerate}

We organize the remaining part of this paper as follows. Section \ref{sec:MFUC} describes the  mathematical formulation of the system's real-time flexibility based on the AGC constrained economic dispatch. Section \ref{sec:CaseSt}  reports the studies of real-time flexibility in IEEE standard systems. Section \ref{sec:conclu} concludes the findings of our paper.  


\section{Mathematical Formulation of  Real-Time Flexibility} \label{sec:MFUC}
{In this section, we deploy a robust optimization based framework to measure the system flexibility of the economic dispatch with the dynamic AGC. In current practice, within every 5-minute interval ED run, the system frequency is adjusted by AGC to its nominal value for every 2-6 seconds. In this paper, we consider the system dynamics within \lred{an} ED run cycle, i.e., 5 minutes. That is, the overall time horizon is set to be 5 minutes. Within the ED cycle, we consider all the AGC cycles as a set of time interval $\Tt$. Both load at each bus and system load disturbance in each time period are assumed to be random and are within \lgreen{an} undetermined variation range. That is,
\begin{align}
& {d_b} \in \mathcal{U}_b(\lambda) = [\bar{d}_{b} - \lambda^{dn}_{b}\hat{d}_b, \bar{d}_b + \lambda^{up}_{b}\hat{d}_b], \ \ \forall b \in \mathcal{B} \nonumber
\end{align}
and
\begin{align}
& \Delta{d}_{t} \in \mathcal{U}_t(\lambda) = [{\Delta\bar{d}}_{t} - \lambda^{dn}_{t}\Delta\hat{d}_t, \Delta\bar{d}_t + \lambda^{up}_{t}\Delta\hat{d}_t], \ \ \forall t \in \Tt, \nonumber
\end{align} 
where $\bar{d}_{b}$ and $\hat{d}_b$ represent the nominal value and maximum deviation of the load at bus $b$, and $\Delta\bar{d}_{t}$ and $\Delta\hat{d}_t$ represent the nominal value and maximum deviation of the system load change at time period $t$. $\lambda_{\{b,t\}}^{\{up, dn\}} \in [0,1]$ represent the scales of the deviation for the corresponding variation range. Unlike the traditional robust optimization model that the size of the uncertainty set is predefined, in our model, we will obtain the largest size of each uncertainty set by deciding the scales $\lambda_{\{b,t\}}^{\{up, dn\}}$, to measure the system flexibility \lred{that} the system can accommodate. 

\lred{Based on the economic dispatch with dynamic AGC model proposed by \cite{Zhang2019AGCED}, we develop a  flexibility measurement model as described in \ref{subsec:FRCTG}.}
\subsection{Formulation}\label{subsec:FRCTG}
\begin{subequations}
\begin{align}
\max \ &  \sum_{\zb{b\in \Bb}}(\hat{d}_b\lambda^{up}_{b} + \hat{d}_b\lambda^{dn}_{b}) + \sum_{t \in \Tt}(\Delta\hat{d}_t\lambda^{up}_{t} + \Delta\hat{d}_t\lambda^{dn}_{t}) \label{eqn-UC:obj}\\
s.t. &  \sum_{n \in {\mathcal{{G}}}}  oc_{n}(p_{n}) + \sum_{n \in \Gg}\sum_{t \in \Tt}\fred{\mbox{CP}_{n}\Delta{f}_{n,t}^{\mbox{\tiny GV}+}} \nonumber \\
	&  \fred{+ \sum_{\zb{n} \in \Gg}\sum_{t \in \Tt} \mbox{CP}_{n}\Delta{f}_{n,t}^{\mbox{\tiny GV}-}} \leq \tau \label{eqn-UC:budget}\\
& p_{n} + reg^{\mbox{\tiny U}}_{n} + sr_{n} \leq \mbox{P}_{n}^{\mbox{\tiny max}}, \forall n \in {\mathcal{{G}}}, \label{eqn-ED:PGUBD}\\
& \mbox{P}_{n}^{\mbox{\tiny min}} \leq p_{n} - reg^{\mbox{\tiny D}}_{n},\forall n \in {\mathcal{{G}}}, \label{eqn-ED:PGLBD}\\
& reg^{\mbox{\tiny U}}_{n} \leq \mbox{REG}^{\mbox{\tiny U}}_{n},\forall n \in {\mathcal{{G}}}, \label{eqn-ED:REGUP}\\
& reg^{\mbox{\tiny D}}_{n} \leq \mbox{REG}^{\mbox{\tiny D}}_{n},\forall n \in {\mathcal{{G}}}, \label{eqn-ED:REGDN} \\
& sr_{n} \leq \mbox{SR}_{n}^{\mbox{\tiny max}},\forall n \in \zb{\mathcal{{G}}}, \label{eqn-ED:SRUBD} \\
& \Delta{p}_{n,t+1}^{\mbox{\tiny M}} = \sum_{i \in \mathcal{G}} (\alpha_{i,n}\Delta{p}_{i,t}^{\mbox{\tiny M}} + \beta_{i,n} \Delta{p}_{i,t}^{\mbox{\tiny GV}})  + \gamma_{n} \Delta{\omega}_{t}\nonumber\\
& + \zeta_{n}\Delta{d}_{t}, \forall n \in {\mathcal{{G}}}, \forall \Delta{d}_{t} \in \mathcal{U}_t,  \forall t \in \Tt,  \label{AGC1}\\
& \Delta{\omega}_{t+1} = \sum_{i \in \mathcal{G}} (\kappa_{i}\Delta{p}_{i,t}^{\mbox{\tiny M}} + \tau_{i}\Delta{p}_{i,t}^{\mbox{\tiny GV}}) + \rho\Delta{\omega}_{t} + \eta\Delta{d}_{t}, \nonumber \\ 
&\forall \Delta{d}_{t} \in \mathcal{U}_t, \forall t \in \Tt, \label{AGC2}\\
& \Delta{p}_{n,t + 1}^{\mbox{\tiny GV}} - \Delta{p}_{n,t}^{\mbox{\tiny GV}} + \fred{\Delta{f}_{n,t + 1}^{\mbox{\tiny GV}+} - \Delta{f}_{n,t+1}^{\mbox{\tiny GV}-}} \nonumber \\
&= K_{n}\Delta{\omega}_{\fred{t+1}},\forall n \in {\mathcal{{G}}}, \forall t \in \Tt \label{AGC3}\\
&  \Delta{p}_{n,t + 1}^{\mbox{\tiny M}} - \Delta{p}_{n,t}^{\mbox{\tiny M}} \leq \mbox{RUR}_{n},\forall n \in {\mathcal{{G}}}, \forall t \in \Tt \label{AGC_up}\\
&  \Delta{p}_{n,t}^{\mbox{\tiny M}} - \Delta{p}_{n,t+1}^{\mbox{\tiny M}} \leq \mbox{RDR}_{n},\forall n \in {\mathcal{{G}}}, \forall t \in \Tt \label{AGC_down}\\
&  -reg^{\mbox{\tiny D}}_{n} \leq  \Delta{p}_{n,t}^{\mbox{\tiny GV}} \leq reg^{\mbox{\tiny U}}_{n}, \forall n \in {\mathcal{{G}}}, \forall t \in \Tt \label{reg}\\
&   \Delta{\omega}^{\mbox{\tiny min}}\leq \Delta{\omega}_{t} \leq \Delta{\omega}^{\mbox{\tiny max}}, \zb{\forall t \in \Tt}\label{frequency}\\
& \sum_{n\in {\mathcal{{G}}}}{p}_{n} - \sum_{b \in \Bb}{d}_{b}=0, \forall {d_b} \in \mathcal{U}_b, \label{eqn-ED-B:Pbalance} \\
&\sum_{n \in \mathcal{{G}}} sr_{n} \geq \mbox{SR}^{\mbox{\tiny min}}, \label{eqn-ED-B: spinning} \\
& \fred{\sum_{n \in \mathcal{{G}}} reg^{\mbox{\tiny U}}_{n} \geq \mbox{REG}^{\mbox{\tiny MinU}}}, \label{eqn-ED-B: Regup} \\
& \fred{\sum_{n \in \mathcal{{G}}} reg^{\mbox{\tiny D}}_{n} \geq \mbox{REG}^{\mbox{\tiny MinD}}},\label{eqn-ED-B: Regdn} \\
& - F_{l} \leq \sum_{b \in \mathcal{B}} \mbox{SF}_{b,l}(\sum_{n \in \mathcal{G}^{b}}p_{n} - {d}_{b}) \leq F_{l}, \forall {d_b} \in \mathcal{U}_b, \label{transmissioncap}\\
&  p_{n},  reg^{\mbox{\tiny U}}_{n},  reg^{\mbox{\tiny D}}_{n}, {sr_{n}}, \zb{\Delta{f}_{n,t}^{\mbox{\tiny GV}+}, \Delta{f}_{n,t}^{\mbox{\tiny GV}-}} \geq 0, \nonumber\\
& \Delta{p}^{\mbox{\tiny GV}}_{n,t},\Delta{p}^{\mbox{\tiny M}}_{n,t},\Delta{\omega}_{t} \ \mbox{free},\\
& \lambda_b^{up}, \lambda_b^{dn},  \lambda_t^{up}, \lambda_t^{dn} \in [0,1], \zb{\forall b \in \Bb, \forall n \in {\mathcal{{G}}}}, \forall t \in \Tt, \nonumber
\end{align}
\end{subequations}
where the objective  function is to maximize the variation range of the uncertainty. Constraints \eqref{eqn-UC:budget} represent the budget constraint, which indicates that the total fuel cost and the penalty cost should not exceed a budget $\tau$. Constraints (\ref{eqn-ED:PGUBD}) and (\ref{eqn-ED:PGLBD}) represent {the generation} limits of {traditional} thermal units that \lgreen{take} account of generation output, regulation and spinning reserves. Constraints (\ref{eqn-ED:REGUP}), (\ref{eqn-ED:REGDN}), and (\ref{eqn-ED:SRUBD}) represent the capacities for providing regulation up, regulation down, and spinning reserve services respectively. Constraints \eqref{AGC1} - \eqref{AGC3} represent AGC dynamic system constraints, i.,e., the transformation of state vectors $\Delta{p}_{n,t}^{\mbox{\tiny M}}$, $\Delta{p}_{n,t}^{\mbox{\tiny GV}}$, $\Delta{\omega}_t$ from $t$ to $t+1$ given any demand disturbance $\Delta{d}_t$, where the matrix components $\alpha, \beta, \gamma, \zeta, \kappa, \tau, \rho, \eta$ can be calculated by numerous methods such as zero-order hold method \cite{Zhang2019AGCED}. Constraints \eqref{AGC_up} and \eqref{AGC_down} \lred{restrict} ramping up and ramping down limits. Constraints \eqref{reg} indicate that the governor generation change should not exceed the regulation service reserved, and constraints \eqref{frequency} restrict the limit of system frequency change.  
In addition, the power balance constraints are described in (\ref{eqn-ED-B:Pbalance}), which should be held for any load realization within the uncertainty set; constraints
(\ref{eqn-ED-B: spinning})-\eqref{eqn-ED-B: Regdn} describe the overall spinning reserve, regulation up, and regulation down requirements respectively, and constraints \eqref{transmissioncap} represent the transmission capacity constraints.

\subsection{Solution Methodology}
First, for notation brevity, we use matrices and vectors to represent constraints and variables, and rewrite the above model in \lgreen{an} abstract compact form (denoted as ACF):
\begin{subequations}
\begin{align}
(\mbox{ACF})\ \  & \max \ \zb{a^T} \lambda \label{primal_abs_obj}\\
\hspace{-0.1cm}s.t.\ \ &  A_1 x \leq b_{1}, \label{primal_cons1} \\
& A_2 x = H_2 d, \ \ \zb{\forall d_b \in \mathcal{U}_b(\lambda), \forall b \in \Bb}, \label{primal_cons6}\\
&A_3 x \leq H_3 d, \ \ \zb{\forall d_b \in \mathcal{U}_b(\lambda), \forall b \in \Bb}, \label{primal_cons7}\\
&\zb{A_{4}y = H_{4} \Delta d, \ \ \Delta d_t \in \mathcal{U}_t(\lambda), \zb{\forall t \in \Tt},} \label{primal_cons2}\\
&\zb{A_{5}y = b_{5},  }\label{primal_cons3}\\
&\zb{A_{6} y \leq b_{6}, }\label{primal_cons4}\\
&\zb{A_7x + A_{8} y \leq b_{7},} \label{primal_cons5}
\end{align}
\end{subequations}
where \zb{$x = (p, reg^{\mbox{\tiny U}}, reg^{\mbox{\tiny D}}, sr)$, $y = (\Delta p^{\mbox{\tiny M}}$, $\Delta{p}^{\mbox{\tiny GV}}$, $\Delta{\omega}$, $\Delta{f}^{\mbox{\tiny GV}+}$, $\Delta{f}^{\mbox{\tiny GV}-})$, $d = (d_1, \cdots,$ $d_b, \cdots)_{b \in \Bb}$, and $\Delta d = (\Delta d_1$, $\cdots$, $\Delta d_t$, $\cdots)_{t \in \Tt}$}; objective \eqref{primal_abs_obj} \lgreen{represents} \eqref{eqn-UC:obj}; constraint \eqref{primal_cons1} represents \eqref{eqn-ED:PGUBD} - \eqref{eqn-ED:SRUBD} and \eqref{eqn-ED-B: spinning} - \eqref{eqn-ED-B: Regdn}; constraint \eqref{primal_cons6} represents \eqref{eqn-ED-B:Pbalance}; constraint \eqref{primal_cons7} represents \eqref{transmissioncap}; constraint \eqref{primal_cons2} represents \eqref{AGC1} and \eqref{AGC2}; constraint \eqref{primal_cons3} represents \eqref{AGC3}; constraint \eqref{primal_cons4} represents \eqref{AGC_up}, \eqref{AGC_down}, and \eqref{frequency}; constraint \eqref{primal_cons5} represents \eqref{eqn-UC:budget} and \eqref{reg}. 

We deploy Benders' decomposition framework to solve the problem. Since the problem is to find the largest deviation of the uncertainty set that the system can accommodate, i.e., the largest value of $\lambda$ without making the constraints \eqref{primal_cons1}-\eqref{primal_cons5} infeasible, therefore, only feasibility cuts are needed. 
\subsubsection{Master Problem and Subproblem}
We first decompose problem ACF into a master problem   \lred{(denoted as MAP)} and a subproblem. 
\begin{subequations}
\begin{align}
(\mbox{MAP})~~ \max \ & \zb{a^T} \lambda \label{primal_MAP_obj}\\
s.t.~~&  g(\lambda) \leq 0 \label{primal_MAP_cons1} \\
& \textbf{0} \leq \lambda \leq \zb{\textbf{1}}
\end{align}
\end{subequations}
Here, the feasibility cuts are represented in \eqref{primal_MAP_cons1}, \zb{and \textbf{0}, $\textbf{1}$ represent vectors with all components 0 and 1 respectively}. 
To generate the feasibility cuts, we first describe the feasibility check problem \lred{ (denoted as FEA)} as follows: 
\begin{subequations}
\begin{align}
(\mbox{FEA}) & \max_{\zb{d_b} \in \mathcal{U}_b, \Delta \zb{d_t} \in \mathcal{U}_t}  \min_{x,y,s}  \ \zb{\textbf{1}}^T s \label{Fea_abs_obj}\\
s.t. \ \  & \fred{A_1} x - s_1\leq b_{1}, \label{Fea_cons1} \\
  &A_2 x +s_{2}^+ - s_{2}^- = H_2 d, \label{Fea_cons6}\\
 &A_3 x  - \zb{s_3} \leq H_3 d, \label{Fea_cons7}\\
 &\zb{A_{4}}y + s_{4}^+ - s_{4}^-  = H_{4} \Delta d, \label{Fea_cons2}\\
 &\fred{\zb{A_{5}}y + s_{5}^+ - s_{5}^- = b_{5}}, \label{Fea_cons3}\\
 &\zb{A_{6}} \zb{y} - s_{6} \leq b_{6}, \label{Fea_cons4}\\
 &A_7x + {\zb{A_{8}}} \zb{y} - s_{7} \leq b_{7}. \label{Fea_cons5}
 \end{align}
 \end{subequations}
If $\lambda$ is feasible, then the optimal value of (FEA) will be $0$.

\subsubsection{Reformulation of Subproblem}
 Now we take the dual of the inner minimization of subproblem (FEA) and combine the dual problem with the outer maximization problem, then we can get the following formulation:
\begin{subequations}
\begin{align}
(\mbox{DFEA}) &\max_{\zb{d_b} \in \mathcal{U}_b, \Delta \zb{d_t} \in \mathcal{U}_t, \mu}  \   \zb{b_1^T \mu_1  + d^T H_2^T \mu_2  + d^T H_3^T \mu_3}   \nonumber \\
& \zb{+  \Delta d^T H_{4}^T \mu_{4} + b_{5}^T \mu_{5} + b_{6}^T \mu_{6} + b_{7}^T \mu_{7}}\\
s.t.\ &  \zb{A_1^T \mu_1 + A_2^T \mu_2 + A_3^T \mu_3 + A_{7}^T \mu_{7}}  \leq 0, \\
& \zb{A_{4}^T}\mu_{4} + \zb{A_{5}^T}\mu_{5} + \zb{A_{6}^T}\mu_{6} + \zb{A_{8}^T}\mu_{7} = 0, \\
& \fred{ -\zb{\textbf{1}} \leq \mu_1, \mu_3, \mu_6, \mu_7 \leq \textbf{0},  } \\ 
& \fred{ -\zb{\textbf{1}} \leq \mu_2, \mu_{4}, \mu_{5} \leq \zb{\textbf{1}}},
\end{align}
\end{subequations}
where $\mu_1$, $\mu_2$, $\mu_3$, $\mu_{4}$, $\mu_{5}$, $\mu_{6}$, $\mu_{7}$ are dual variables for constraints \eqref{Fea_cons1} - \eqref{Fea_cons5} respectively.

\fred{In the above formulation (DFEA), we have  bilinear terms \zb{$d^T H_2^T \mu_2$, $d^T H_3^T \mu_3$, and $\Delta d^T H_{4}^T \mu_{4}$}. We will deal with \zb{$d^T H_3^T \mu_3$} first. \zb{Let $N_i$ represent the dimension of $\mu_i$ and $\lambda^{*}$ represent} the optimal solution of problem (MAP). Based on the property of $d$, \zb{$\forall b \in \Bb$,} we can rewrite it as \zb{$d_b = \bar{d}_b + z^{+}_b \lambda^{*,up}_b \hat{d}_b -  z^{-}_b \lambda^{*,dn}_b \hat{d}_b$}. Here we introduce two binary variables $\zb{z^{+}_b}$ and $z^{-}_b$ to indicate the deviation direction. Note that the variables  $\lambda^{*,up}$ and $\lambda^{*,dn}$ have been fixed for DFEA problem.  Therefore, we can replace the bilinear term $d^T H_3^T \mu_3$ as follows:}
\begin{subequations}
\begin{align}
& d^T H_3^T \mu_3 = \sum_{b\in \Bb}\sum_{i=1}^{N_3} d_b \lgreen{H_{3,i,b}}\mu_{3,i}\nonumber\\
& = \sum_{b \in \Bb}\sum_{i=1}^{N_3}(\bar{d}_b \lgreen{H_{3,i,b}} \mu_{3,i} + \lambda^{*,up}_b\hat{d}_b \lgreen{H_{3,i,b}} \zb{z^{+}_b}\mu_{3,i} \nonumber\\
&- \lambda^{*,dn}_b\hat{d}_b \lgreen{H_{3,i,b}} \zb{z^{-}_b}{\mu}_{3,i}) \label{Refm_dmu}\\
& \zb{z^{+}_b} + \zb{z^{-}_b} = 1, z_b^{+}, z_b^{-} \in \{0,1\} \label{Refm_unique}
\end{align}
\end{subequations}
\fred{
In (\ref{Refm_dmu}), we have two bilinear items $\zb{z^{+}_b}\mu_{3,i}$ and $\zb{z^{-}_b}\mu_{3,i}$, which can be further linearized by introducing auxiliary variables $\mu_{3,b,i}^{+}$ and $\mu_{3,b,i}^{-}$. 
By following the approach indicated in \cite{Fan2014MinMaxRegret}, we have the following reformulation which is equivalent to (\ref{Refm_dmu}):}
\zb{
\begin{subequations}
\begin{align}
&  d^T H_3^T \mu_3 = \sum_{b \in \Bb}\sum_{i=1}^{N_3}(\bar{d}_b \lgreen{H_{3,i,b}}\mu_{3,i} + \nonumber\\
&\lambda^{*,up}_b\hat{d}_b \lgreen{H_{3,i,b}} \mu_{3,b,i}^{+} - \lambda^{*,dn}_b\hat{d}_b \lgreen{H_{3,i,b}} {\mu}_{3,b,i}^{-}) \label{Refm_dmuR}\\
&-\zb{z^{+}_b} \leq \mu_{3,b,i}^{+}, \ \mu_{3,i} \leq \mu_{3,b,i}^{+} \leq 1 - \zb{z^{+}_b} + \mu_{3,i},  \\
&-\zb{z^{-}_b} \leq \mu_{3,b,i}^{-},\ \mu_{3,i} \leq \mu_{3,b,i}^{-} \leq 1 - \zb{z^{-}_b} + \mu_{3,i}, \\
& z_b^+ + z_b^- = 1, z_b^{+}, z_b^{-} \in \{0,1\},  \\
&-1 \leq \mu_{3,b,i}^{+}, \mu_{3,b,i}^{-} \leq 0, \forall b \in \Bb, \forall i = 1, \cdots, N_3. \label{last}
\end{align}
\end{subequations}
}
For \zb{$d^T H_2^T \mu_2$ and $\Delta d^T H_{4}^T \mu_{4}$}, we can use a similar approach. \zb{We will use $d^T H_2^T \mu_2$ as an example and the other follows the same approach. Since $\mu_2 \in [-\textbf{1},\textbf{1}]$, we can replace it with $\mu_{2}^n - \mu_{2}^p$ and $- \textbf{1} \leq \mu_{2}^n \leq \textbf{0}$ and $-\textbf{1}\leq \mu_{2}^p \leq \lgreen{\textbf{0}} $. Then $d^T H_2^T\mu_2 = d^T H_2^T\mu_{2}^n - d^T H_2^T\mu_{2}^p$. For $d^T H_2^T\mu_{2}^n$ and $d^T H_2^T\mu_{2}^p$}, we will follow the same procedure of \eqref{Refm_dmuR}-\eqref{last} to linearize them. Therefore, we can reformulate the (DEFA) problem as follows:
\zb{
\begin{subequations}
\begin{align}
(\mbox{RDFEA}) \ &\max_{d_b, \Delta d_t, \mu}  \   b_1^T \mu_1   +\sum_{b \in \Bb}\sum_{i=1}^{N_2} \{ \bar{d}_b \lgreen{ H_{2,i,b}} (\mu_{2,i}^n-\mu_{2,i}^p)+ \nonumber\\
& \hspace{-1.6cm}\lambda^{*,up}_b\hat{d}_b\lgreen{ H_{2,i,b}} (\mu_{2,b,i}^{n,+}-\mu_{2,b,i}^{p,+}) - \lambda^{*,dn}_b\hat{d}_b \lgreen{ H_{2,i,b}} ({\mu}_{2,b,i}^{n,-}-{\mu}_{2,b,i}^{p,-})\} \nonumber\\
&\hspace{-1.6cm}  + \sum_{b \in \Bb}\sum_{i=1}^{N_3}(\bar{d}_b \lgreen{ H_{3,i,b}} \mu_{3,i} + \lambda^{*,up}_b\hat{d}_b \lgreen{ H_{3,i,b}}\mu_{3,b,i}^{+} \nonumber\\
&\hspace{-1.6cm} - \lambda^{*,dn}_b\hat{d}_b \lgreen{ H_{3,i,b}} {\mu}_{3,b,i}^{-}) + \sum_{t \in \Tt}\sum_{i=1}^{N_4} \{ \Delta \bar{d}_t \lgreen{ H_{4,i,t}}(\mu_{4,i}^n-\mu_{4,i}^p)  \nonumber\\
& \hspace{-1.6cm}+\lambda^{*,up}_t \Delta\hat{d}_t \lgreen{ H_{4,i,t}}(\mu_{4,t,i}^{n,+}- \mu_{4,t,i}^{p,+}) \nonumber\\
& \hspace{-1.6cm}- \lambda^{*,dn}_b \lgreen{\Delta\hat{d}_t} \lgreen{ H_{4,i,t}} ({\mu}_{4,t,i}^{n,-}-{\mu}_{4,t,i}^{p,-})\} + b_{5} \mu_{5}  + b_{6} \mu_{6} + b_{7} \mu_{7} \nonumber\\
& \hspace{-1cm}  s.t. \ \ \ \lgreen {A_1^{T} \mu_1 + A_2^{T} \mu_2 + A_3^{T} \mu_3 + A_{7}^{T} \mu_{7}  \leq 0}, \\
& \hspace{-0.8 cm} \lgreen{ \zb{A_{4}^{T}}\mu_{4} + \zb{A_{5}^{T}}\mu_{5} + \zb{A_{6}^{T}}\mu_{6} + \zb{A_{8}^{T}}\mu_{7} = 0},\\
&\hspace{-0.8cm}  -\zb{z^{+}_b} \leq \mu_{2,b,i}^{n,+}, \ \mu_{2,i}^n \leq \mu_{2,b,i}^{n,+} \leq 1 - \zb{z^{+}_b} + \mu_{2,i}^n,\\
&\hspace{-0.8cm}  -\zb{z^{-}_b} \leq \mu_{2,b,i}^{n,-},\  \mu_{2,i}^n \leq \mu_{2,b,i}^{n,-} \leq 1 - \zb{z^{-}_b} + \mu_{2,i}^n,\\
&\hspace{-0.8cm}  -\zb{z^{+}_b} \leq \mu_{2,b,i}^{p,+},\  \mu_{2,i}^p \leq \mu_{2,b,i}^{p,+} \leq 1 - \zb{z^{+}_b} + \mu_{2,i}^p,\\
&\hspace{-0.8cm}  -\zb{z^{-}_b} \leq \mu_{2,b,i}^{p,-}, \ \mu_{2,i}^p \leq \mu_{2,b,i}^{p,-} \leq 1 - \zb{z^{-}_b} + \mu_{2,i}^p,\\
&\hspace{-0.8cm} \forall b \in \Bb, \forall i = 1, \cdots, N_2, \nonumber\\
&\hspace{-0.8cm}  -\zb{z^{+}_b} \leq \mu_{3,b,i}^{+},\  \mu_{3,i} \leq \mu_{3,b,i}^{+} \leq 1 - \zb{z^{+}_b} + \mu_{3,i},\\
&\hspace{-0.8cm}  -\zb{z^{-}_b} \leq \mu_{3,b,i}^{-},\  \mu_{3,i} \leq \mu_{3,b,i}^{-} \leq 1 - \zb{z^{-}_b} + \mu_{3,i},\\
&\hspace{-0.8cm} \forall b \in \Bb, \forall i = 1, \cdots, N_3, \nonumber\\
&\hspace{-0.8cm}  -\zb{z^{+}_t} \leq \mu_{4,t,i}^{n,+},\  \mu_{4,i}^n \leq \mu_{4,t,i}^{n,+} \leq 1 - \zb{z^{+}_t} + \mu_{4,i}^n,\\
&\hspace{-0.8cm}  -\zb{z^{-}_t} \leq \mu_{4,t,i}^{n,-},\  \mu_{4,i}^n \leq \mu_{4,t,i}^{n,-} \leq 1 - \zb{z^{-}_t} + \mu_{4,i}^n,\\
&\hspace{-0.8cm}  -\zb{z^{+}_t} \leq \mu_{4,t,i}^{p,+},\  \mu_{4,i}^p \leq \mu_{4,t,i}^{p,+} \leq 1 - \zb{z^{+}_t} + \mu_{4,i}^p,\\
&\hspace{-0.8cm}  -\zb{z^{-}_t} \leq \mu_{4,t,i}^{p,-},\ \mu_{4,i}^p \leq \mu_{4,t,i}^{p,-} \leq 1 - \zb{z^{-}_t} + \mu_{4,i}^p,\\
&\hspace{-0.8cm} \forall t \in \Tt, \forall i = 1, \cdots, N_4, \nonumber\\
&\hspace{-0.8cm}   z_b^+ + z_b^- = 1, \ \ z_t^+ + z_t^- = 1, \forall b \in \Bb, \forall t \in \Tt, \\
&\hspace{-0.8cm}  -\textbf{1} \leq\mu_{2}^{n,+}, \mu_{2}^{n,-}, \mu_{2}^{p,+}, \mu_{2}^{p,-}  \leq \textbf{0}, \\
&\hspace{-0.8cm} -\textbf{1}\leq  \mu_{3}^{+}, \mu_{3}^{-} \leq \textbf{0}, \nonumber\\
&\hspace{-0.8cm}  -\textbf{1} \leq\mu_{4}^{n,+}, \mu_{4}^{n,-} , \mu_{4}^{p,+},\mu_{4}^{p,-} \leq \textbf{0}, \\
&\hspace{-0.8cm}   -\textbf{1} \leq \mu_1, \mu_3, \mu_6, \mu_7 \leq \textbf{0},  \\ 
&\hspace{-0.8cm}   -\textbf{1} \leq \mu_2, \mu_{4}, \mu_{5} \leq \textbf{1}. 
\end{align}
\end{subequations}
}
If the objective value $\eta^{*}$ of subproblem (RDFEA) is greater than $0$, we can obtain the following feasibility cut:
%
%
\zb{
\begin{align}
& \eta(\lambda) = \nonumber \\
&\  \sum_{b \in \Bb}(\sum_{i=1}^{N_2}\hat{d}_b \lgreen{H_{2,i,b}}(\mu_{2,b,i}^{n+,*}-\mu_{2,b,i}^{p+,*}) +  \sum_{i=1}^{N_3}\hat{d}_b\lgreen{H_{3,i,b}}\mu_{3,b,i}^{+,*})\lambda^{up}_b   \nonumber\\
& -  \sum_{b \in \Bb}(\sum_{i=1}^{N_2}\hat{d}_b \lgreen{H_{2,i,b}}(\mu_{2,b,i}^{n-,*}-\mu_{2,b,i}^{p-,*}) +  \sum_{i=1}^{N_3}\hat{d}_b \lgreen{H_{3,i,b}}\mu_{3,b,i}^{-,*})\lambda^{dn}_b   \nonumber\\
&+\sum_{t \in \Tt}\sum_{i=1}^{N_4} \Delta\hat{d}_t \lgreen{H_{4,i,t}}(\mu_{4,t,i}^{n+,*}-\mu_{4,t,i}^{p+,*})\lambda^{up}_t\nonumber\\
& - \sum_{t \in \Tt}\sum_{i=1}^{N_4} \lgreen{\Delta\hat{d}_t}\lgreen{H_{4,i,t}} ({\mu}_{4,t,i}^{n-,*}-{\mu}_{4,t,i}^{p-,*})\lambda^{dn}_t \nonumber\\
& +b_1^T \mu_1^{*}   +\sum_{b \in \Bb}\sum_{i=1}^{N_2} \bar{d}_b \lgreen{H_{2,i,b}}(\mu_{2,i}^{n*}-\mu_{2,i}^{p*})\nonumber\\
& + \sum_{b \in \Bb}\sum_{i=1}^{N_3}\bar{d}_b \lgreen{H_{3,i,b}}\mu_{3,i}^{*} + \sum_{t \in \Tt}\sum_{i=1}^{N_4}  \Delta\bar{d}_t\lgreen{H_{4,i,t}}(\mu_{4,i}^{n*}-\mu_{4,i}^{p*})\nonumber\\
& + b_{5} \mu_{5}^{*}  + b_{6} \mu_{6}^{*} + b_{7} \mu_{7}^{*} \leq 0 \label{feaCUT}
\end{align}}
\subsubsection{Algorithm Framework}
We summarize our Benders' decomposition based framework as follows:
\begin{enumerate}
	\item Solve the master problem (MAP) and obtain the optimal solution $\lambda^*$.
	\item Test the feasibility of the subproblem by solving dual reformulated DFEA problem with $\lambda^*$.
	\item If the optimal value $\eta^{*}$ of DFEA $ > 0$, then update the master problem by adding feasibility cut \eqref{feaCUT} and go to Step 2. Otherwise, terminate and return the optimal solution and objective value of master problem.
\end{enumerate}

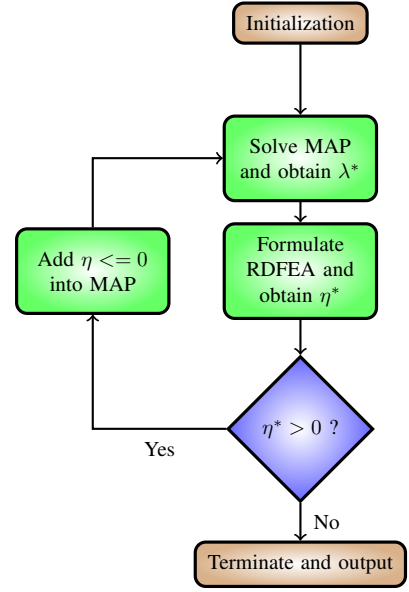
\begin{figure}[tph]
	\begin{center}
		
		\usetikzlibrary{shapes,arrows}
		
		\tikzstyle{decision} = [diamond, draw,,outer color=blue!60,inner color=white,very thick, fill=white!20, text
		width=6em, text badly centered, node distance=4cm, inner sep=0pt]
		\tikzstyle{smalldecision} = [diamond, draw, fill=white!20, text
		width=4em, text badly centered, node distance=4cm, inner sep=0pt]
		\tikzstyle{block} = [rectangle, draw,outer color=brown!60,inner color=white,very thick, fill=white!20, rounded
		corners, minimum size=1em] \tikzstyle{bigblock} = [rectangle, draw,
		fill=white!20, rounded corners, text width=14em, minimum height=3em]
		\tikzstyle{midblock} = [rectangle, draw,,outer color= green!60,inner color=white,very thick, fill=white!20, rounded
		corners, text width=6em, minimum height=4em] \tikzstyle{line} =
		[draw, ->, thick, color=black] \tikzstyle{cloud} = [draw,
		ellipse,fill=black, minimum size=1em]
		\tikzstyle{point}=[circle,inner sep=0pt,minimum size=2pt,fill=black]
		\tikzstyle{every node}=[font=\large]
		\tikzset{
			every node/.style={scale=0.8}
		}
		
		\begin{tikzpicture}[scale=1.2, node distance = 2cm, inner sep=6pt, text badly centered,auto]
		\node [block] (init) at (2.5,10.5) {{Initialization}};
		\node [midblock] (Smaster) at (2.5,9) {Solve $\mbox{MAP}$ and {obtain $\lambda^*$}};
		\node [midblock] (Gsub) at (2.5,7.75) {Formulate $\mbox{RDFEA}$ and {obtain $\eta^{*}$}};
		\node [decision] (Cut) at (2.5, 6) {$\eta^{*} > 0$ ?};
		\node [midblock] (Update) at (0.2, 7.75){Add $\eta <= 0$ into $\mbox{MAP}$};
		\node [block] (stop) at (2.5, 4.5){Terminate and output};
		\path [line] (init) -- (Smaster);
		\path [line] (Smaster) -- (Gsub);
		\path [line] (Gsub) -- (Cut);
		\path [line] (Update) |- (Smaster);
		\path [line] (Cut) -| node [near start] {Yes} (Update);
		\path [line] (Cut) -- node {No}(stop);
		
		\end{tikzpicture}
	\end{center}
	\caption{{{Flowchart} of Decomposition Algorithm}}\label{ALG:BendersChart}
\end{figure}

\section{Case Study}\label{sec:CaseSt}
In this section, we test the performance of the proposed real-time flexibility management model with the IEEE 118-bus system (online \cite{IIT_Systems}). These cases are tested by using Julia \cite{lubin2015computing} and CPLEX 12.8 \cite{IBMCPLEX} on Intel Xeon Silver 4216 CPU and 128 G memory.

To get the insights \lgreen{on} how the system operators can improve their system flexibility by using limited resources, we investigate on which factors will have impacts on the system flexibility. More specifically, we study key factors including the budget of operational cost, ramping capability, and transmission line capacity. In order to obtain more information on the system's capability in handling uncertainty, we analyze seven different system flexibility indicators: \lred{total flexibility (TF), economic dispatch flexibility (EDF), AGC flexibility (AGCF), economic dispatch upward flexibility (EDUPF), economic dispatch downward flexibility (EDDNF), AGC upward flexibility (AGCUPF), and AGC downward flexibility (AGCDNF)} based on the simulation results. The definitions of these flexibility indicators are shown in (\ref{eqn-TF}) - (\ref{eqn-AGCDNF}).
\begin{subequations}
  \begin{align}
 & \mbox{TF}  = \sum_{b}(\hat{d}_b\lambda^{up}_{b} + \hat{d}_b\lambda^{dn}_{b}) + \sum_{t \in \Tt}(\Delta\hat{d}_t\lambda^{up}_{t} + \Delta\hat{d}_t\lambda^{dn}_{t}) \label{eqn-TF}\\
 & \mbox{EDF}  = \sum_{b}(\hat{d}_b\lambda^{up}_{b} + \hat{d}_b\lambda^{dn}_{b})  \label{eqn-EDF}\\
 & \mbox{AGCF}  =  \sum_{t \in \Tt}(\Delta\hat{d}_t\lambda^{up}_{t} + \Delta\hat{d}_t\lambda^{dn}_{t}) \label{eqn-AGCF}\\
 & \mbox{EDUPF}  = \sum_{b}(\hat{d}_b\lambda^{up}_{b})  \label{eqn-EDUPF}\\
 & \mbox{EDDNF}  =  \sum_{b}(\hat{d}_b\lambda^{dn}_{b})  \label{eqn-EDDNF}\\
 & \mbox{AGCUPF}  =  \sum_{t \in \Tt}(\Delta\hat{d}_t\lambda^{up}_{t} ) \label{eqn-AGCUPF}\\
 & \mbox{AGCDNF}  =  \sum_{t \in \Tt}(\Delta\hat{d}_t\lambda^{dn}_{t}) \label{eqn-AGCDNF}
  \end{align}
  \end{subequations}
The TF can gauge both static and dynamic system capability to handle the uncertainty caused by forecasting errors of load and renewable energy resources. The EDF reflects the static capability of the system. The AGCF indicates the capability of the system to respond to the disturbance in real time. The EDUPF and EDDNF can explain up and down static flexibility of the system.  The \lred{AGCUPF} and AGCDNF can describe the up and down dynamic flexibility of the system, respectively. 

In the nominal case, there are 30 units online and the system load is  \lgreen{around} $86.3\%$ of the total capacity of online units. We assume the uncertain range of the net load at each bus is $[-15\%, 15\%]$ and the system load disturbance is within the range $[-1\%, 1\%]$.

\subsection{Impact of Operational Cost}
In this subsection, we first investigate the impact of the operational cost budget by increasing it gradually (denoted from B0 to B22 in Column 1, Table \ref{tab:SFUDB}).  We set the budget of the base case (denoted as B0) as the optimal operational cost of the economic dispatch  under the nominal foretasted net load without the dynamic AGC. Then we incrementally scale up the budget by using the scale factor (SF) as shown in Column 2 in Table \ref{tab:SFUDB}.  For example, the budget of B22 is four times of the budget of B0.  All other system parameters are the same in these 23 cases. We report the simulation results of the proposed flexibility indicators under different budgets in Table \ref{tab:SFUDB}.
 
\begin{table}[h]
	\huge
	\centering
	\caption{{System Flexibility under Different Budgets} }
	\label{tab:SFUDB}
	\begin{adjustbox}{max width=0.48\textwidth}
		\begin{tabular}{|c|c|c|c|c|c|c|c|c|}
			\hline
	Budget&	SF & TF      & EDF     & AGCF  & EDUPF  & EDDNF  & AGCUPF & AGCDNF \\ \hline
	B0	&	1.00        & 115.298 & 115.298 & 0.000 & 47.626 & 67.672 & 0.000  & 0.000  \\ \hline
	B1	&	1.01        & 120.678 & 120.678 & 0.000 & 53.006 & 67.672 & 0.000  & 0.000  \\ \hline
	B2	&	1.02        & 125.532 & 125.532 & 0.000 & 57.861 & 67.672 & 0.000  & 0.000  \\ \hline
	B3	&	1.03        & 129.411 & 129.411 & 0.000 & 61.739 & 67.672 & 0.000  & 0.000  \\ \hline
	B4	&	1.04        & 129.520 & 129.476 & 0.044 & 61.804 & 67.672 & 0.026  & 0.018  \\ \hline
	B5	&	1.05        & 129.550 & 129.476 & 0.074 & 61.804 & 67.672 & 0.048  & 0.026  \\ \hline
	B6	&	1.06        & 129.579 & 129.476 & 0.103 & 61.804 & 67.672 & 0.078  & 0.026  \\ \hline
	B7	&	1.07        & 129.609 & 129.476 & 0.133 & 61.804 & 67.672 & 0.059  & 0.074  \\ \hline
	B8	&	1.08        & 129.638 & 129.476 & 0.162 & 61.804 & 67.672 & 0.086  & 0.076  \\ \hline
	B9	&	1.09        & 129.668 & 129.476 & 0.192 & 61.804 & 67.672 & 0.093  & 0.099  \\ \hline
	B10	&	1.10        & 129.697 & 129.476 & 0.221 & 61.804 & 67.672 & 0.109  & 0.112  \\ \hline
	B11	&	1.20        & 129.993 & 129.476 & 0.517 & 61.804 & 67.672 & 0.290  & 0.227  \\ \hline
	B12	&	1.30        & 130.288 & 129.476 & 0.812 & 61.804 & 67.672 & 0.400  & 0.412  \\ \hline
	B13	&	1.40       & 130.583 & 129.476 & 1.107 & 61.804 & 67.672 & 0.578  & 0.529  \\ \hline
	B14	&	1.50        & 130.878 & 129.476 & 1.402 & 61.804 & 67.672 & 0.730  & 0.672  \\ \hline
	B15	&	1.60        & 131.174 & 129.476 & 1.698 & 61.804 & 67.672 & 0.682  & 1.015  \\ \hline
	B16	&	1.70        & 131.469 & 129.476 & 1.993 & 61.804 & 67.672 & 0.844  & 1.149  \\ \hline
	B17	&	1.80        & 131.764 & 129.476 & 2.288 & 61.804 & 67.672 & 1.133  & 1.156  \\ \hline
	B18	&	1.90        & 131.976 & 129.476 & 2.500 & 61.804 & 67.672 & 1.250  & 1.250  \\ \hline
	B19	&	2.00        & 131.976 & 129.476 & 2.500 & 61.804 & 67.672 & 1.250  & 1.250  \\ \hline
	B20	&	2.50        & 131.976 & 129.476 & 2.500 & 61.804 & 67.672 & 1.250  & 1.250  \\ \hline
	B21	&	3.00        & 131.976 & 129.476 & 2.500 & 61.804 & 67.672 & 1.250  & 1.250  \\ \hline
	B22	&	4.00        & 131.976 & 129.476 & 2.500 & 61.804 & 67.672 & 1.250  & 1.250  \\ \hline
		\end{tabular}
	\end{adjustbox}
\end{table}

From Columns 1-5 in Table \ref{tab:SFUDB}, we can observe that three flexibility indicators (i.e., TF, EDF, and \lred{AGCF}) are non-decreasing as the budget increases. 
On one hand, for EDF, when the budget is small and gradually increases by $1\%$ from B0 to B3, the EDF increases by around $4\%$ for each level of the budget. Once \lgreen{the} scale factor (SF) reaches $1.04$, EDF remains as a constant number. On the other hand, when the budget remains small, i.e., B0-B3, AGCF remains 0. But when the budget becomes larger, AGCF increases and becomes a constant \lgreen{when a relatively large budget (i.e., B18 and above) is available.} That is because that when the budget is small, EDF will be allocated to all operational cost budget as it has more weight than AGCF in the objective function. But when the system has more budget for operational cost,  which is sufficient for the operational cost of the economic dispatch, it will start to allocate budget to the dynamic AGC part. Therefore, with more budget available, AGCF will incur more operational \lgreen{costs}. But when the total budget is sufficient for both ED and AGC costs, both EDF and AGCF will not increase. In addition, the increment of AGCF will result in the increment of TF. We show the trends of indicators TF, EDF, and AGCF in Fig. \ref{FIG:sfB}.
\begin{figure}[tph]
\begin{center}
\input{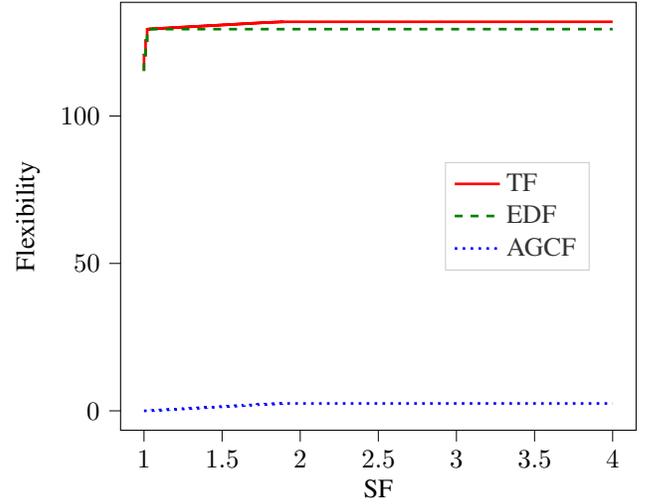}
\end{center}
\caption{Trends of TF, EDF, and AGCF }\label{FIG:sfB}
\end{figure}

We further analyze the trend of EDF by looking into its two components: EDUPF and EDDNF, and we show the \zb{trends} of these three factors in Fig. \ref{FIG:edfB}.  From the \lgreen{figure,} we can observe that EDDNF does not change as the budget changes but \lred{EDUPF} increases and becomes stable \lgreen{when  a higher budget is available}. That is because the operational cost only includes the fuel cost of thermal units, which means that it can reflect the cost caused by increasing the thermal \lgreen{generators'} output. In other words, the economic dispatch \lgreen{upward} flexibility can be restricted by the system budget when the system experiences a high load and a low renewable energy output. Therefore, the budget of the operational cost has a significant impact on EDUPF, especially when the budget is close to the cost of the economic dispatch only without dynamic AGC. On the contrary,  since there is no cost when generators reduce their output, the budget has little impact on the economic dispatch \lgreen{downward} flexibility.  
\begin{figure}[tph]
	\begin{center}
	\input{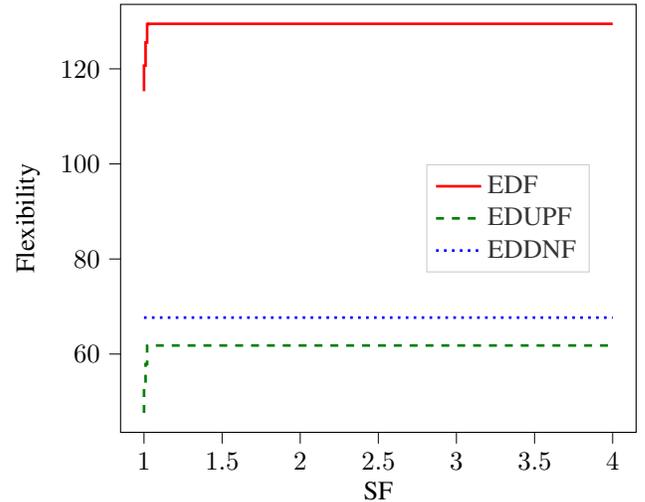}
\end{center}
\caption{Trends of EDF, EDUPF, and EDDNF }\label{FIG:edfB}
\end{figure}

In addition, we show the trends of AGCF, AGCUPF, and AGCDNF in Fig. \ref{FIG:agcB}. We can observe that the AGCF is monotonic non-decreasing as \lgreen{the} budget increases.  Similarly, AGCUPF and AGCDNF both increase (not monotonic) as \lgreen{the} budget increases.  Once the scale factor of the budget reaches B18, all flexibility indicators remain as constants. 

\begin{figure}[tph]
	\begin{center}
\input{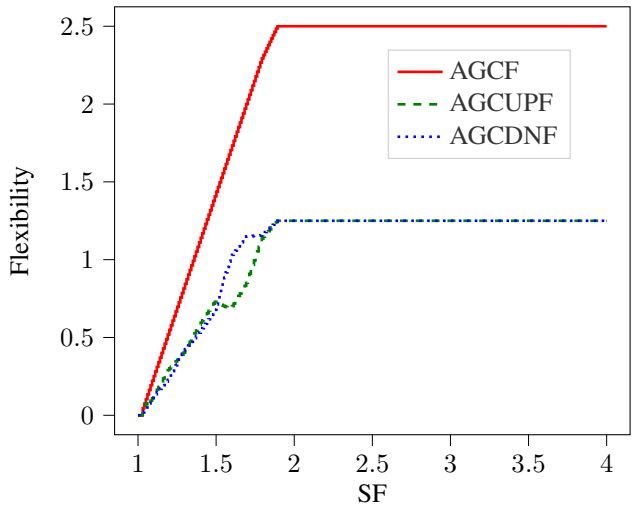}
\end{center}
\caption{Trends of AGCF, AGCUPF, and AGCDNF }\label{FIG:agcB}
\end{figure}

\subsection{Impacts of Ramping Capability}
In this subsection, we report simulation results by increasing the ramping capability of generators under three different budget settings (i.e., B0, B14, and B19).  Then, we test different ramping capability settings of all generators in the system to check how flexibility indicators change with ramping capability under a fixed budget.   We report  flexibility indicators (i.e., TF, EDF, and AGCF) in Table \ref{tab:SFUDRC} and show the TF curves in Fig. \ref{FIG:rampTF}.

Fig. \ref{FIG:rampTF} shows that the increasing ramping capability can improve both TF and EDF of the system, but it has \lgreen{a} limited effect on enhancing AGCF. Under both budgets (B0 and B14), AGCF \zb{increases} as the SF increases, while under budget B19, the AGCF remains the same with the change of ramping capability. Therefore, we can observe that the ramping capacity will \lred{influence} the economic dispatch significantly but \lred{influence} the dynamic AGC \lred{minimally}.
\begin{table}[h]
	\huge
	\centering
	\caption{{System Flexibility under Different Ramping Capabilities} }
	\label{tab:SFUDRC}
	\begin{adjustbox}{max width=0.48\textwidth}
		\begin{tabular}{|c|c|c|c|c|c|c|c|c|c|}
			\hline
			\multirow{2}{*}{SF} & \multicolumn{3}{c|}{B0}   & \multicolumn{3}{c|}{B14}  & \multicolumn{3}{c|}{B19}  \\ \cline{2-10} 
			& TF      & EDF     & AGCF  & TF      & EDF     & AGCF  & TF      & EDF     & AGCF  \\ \hline
			1                   & 115.298 & 115.298 & 0.000 & 130.878 & 129.476 & 1.402 & 131.976 & 129.476 & 2.500 \\ \hline
			1.01                & 115.431 & 115.431 & 0.000 & 130.904 & 129.501 & 1.403 & 132.001 & 129.501 & 2.500 \\ \hline
			1.02                & 115.564 & 115.563 & 0.000 & 130.930 & 129.526 & 1.404 & 132.026 & 129.526 & 2.500 \\ \hline
			1.03                & 115.697 & 115.696 & 0.000 & 130.956 & 129.551 & 1.405 & 132.051 & 129.551 & 2.500 \\ \hline
			1.04                & 115.830 & 115.829 & 0.001 & 130.982 & 129.576 & 1.406 & 132.076 & 129.576 & 2.500 \\ \hline
			1.05                & 115.963 & 115.962 & 0.001 & 131.007 & 129.601 & 1.406 & 132.101 & 129.601 & 2.500 \\ \hline
			1.06                & 116.096 & 116.095 & 0.001 & 131.033 & 129.626 & 1.407 & 132.126 & 129.626 & 2.500 \\ \hline
			1.07                & 116.229 & 116.227 & 0.001 & 131.059 & 129.651 & 1.408 & 132.151 & 129.651 & 2.500 \\ \hline
			1.08                & 116.362 & 116.360 & 0.001 & 131.085 & 129.676 & 1.409 & 132.176 & 129.676 & 2.500 \\ \hline
			1.09                & 116.494 & 116.493 & 0.001 & 131.111 & 129.701 & 1.410 & 132.201 & 129.701 & 2.500 \\ \hline
			1.1                 & 116.623 & 116.622 & 0.002 & 131.136 & 129.726 & 1.410 & 132.226 & 129.726 & 2.500 \\ \hline
			1.2                 & 117.906 & 117.902 & 0.003 & 131.394 & 129.976 & 1.418 & 132.476 & 129.976 & 2.500 \\ \hline
			1.3                 & 119.152 & 119.147 & 0.005 & 131.652 & 130.226 & 1.426 & 132.726 & 130.226 & 2.500 \\ \hline
			1.4                 & 120.367 & 120.361 & 0.006 & 131.909 & 130.476 & 1.433 & 132.976 & 130.476 & 2.500 \\ \hline
			1.5                 & 121.573 & 121.566 & 0.006 & 132.166 & 130.726 & 1.440 & 133.226 & 130.726 & 2.500 \\ \hline
			1.6                 & 122.746 & 122.739 & 0.007 & 132.423 & 130.976 & 1.447 & 133.476 & 130.976 & 2.500 \\ \hline
			1.7                 & 123.864 & 123.857 & 0.007 & 132.680 & 131.226 & 1.454 & 133.726 & 131.226 & 2.500 \\ \hline
			1.8                 & 124.962 & 124.955 & 0.008 & 132.936 & 131.476 & 1.460 & 133.976 & 131.476 & 2.500 \\ \hline
			1.9                 & 126.054 & 126.046 & 0.008 & 133.193 & 131.726 & 1.467 & 134.226 & 131.726 & 2.500 \\ \hline
			2                   & 127.121 & 127.113 & 0.008 & 133.449 & 131.976 & 1.473 & 134.476 & 131.976 & 2.500 \\ \hline
			2.5                 & 131.161 & 131.151 & 0.011 & 133.775 & 132.259 & 1.515 & 134.760 & 132.259 & 2.500 \\ \hline
			3                   & 132.315 & 132.259 & 0.056 & 133.811 & 132.259 & 1.552 & 134.760 & 132.259 & 2.500 \\ \hline
			4                   & 132.393 & 132.259 & 0.133 & 133.874 & 132.259 & 1.615 & 134.760 & 132.259 & 2.500 \\ \hline
		\end{tabular}
	\end{adjustbox}
\end{table}

\begin{figure}[tph]
	\begin{center}
		\input{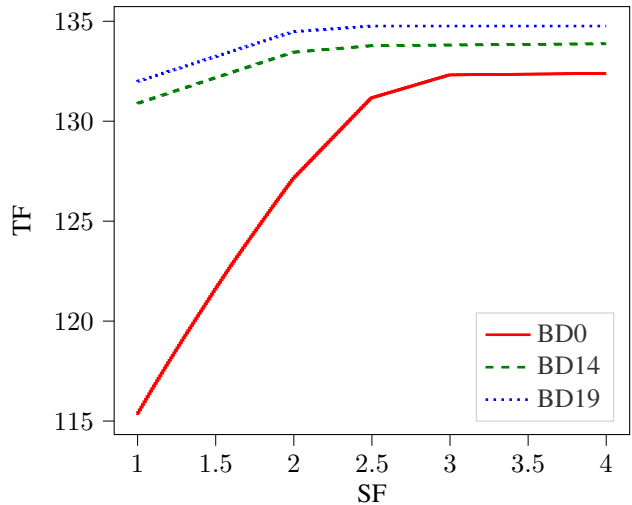}
	\end{center}
	\caption{Trends of TF under Different Ramping Capabilities }\label{FIG:rampTF}
\end{figure}

\subsection{Impacts of Transmission Line Capacity}
In this subsection, we study the impact of transmission line capacity on the system flexibility. Similar \lgreen{to} the setup of ramping capability experiments,  we study the system flexibility under three different budgets (i.e., B0, B14, and B19). Then we test different scenarios of all the transmission line capacity in Table \ref{tab:SFUDTLC}  by changing the scale factor.

We can observe that under the budget B0 and B19, AGCF remains unchanged, and under budget B14, it \lgreen{changes} slightly. \zb{This} also indicates that improving the transmission capacity \zb{has a} marginal impact on the AGCF. On the contrary, Fig. \ref{FIG:lineTF} demonstrates that a larger  transmission line capacity can improve TF for a fixed budget until TF reaches its maximum value.    
\begin{table}[h]
	\huge
	\centering
	\caption{{System Flexibility under Different Transmission Line Capabilities} }
	\label{tab:SFUDTLC}
	\begin{adjustbox}{max width=0.48\textwidth}

	\begin{tabular}{|c|c|c|c|c|c|c|c|c|c|}
		\hline
		\multirow{2}{*}{SF} & \multicolumn{3}{c|}{B0}   & \multicolumn{3}{c|}{B14}  & \multicolumn{3}{c|}{B19}  \\ \cline{2-10} 
		& TF      & EDF     & AGCF  & TF      & EDF     & AGCF  & TF      & EDF     & AGCF  \\ \hline
		1                   & 115.298 & 115.298 & 0.000 & 130.878 & 129.476 & 1.402 & 131.976 & 129.476 & 2.500 \\ \hline
		1.01                & 115.560 & 115.560 & 0.000 & 131.076 & 129.673 & 1.402 & 132.173 & 129.673 & 2.500 \\ \hline
		1.02                & 115.823 & 115.823 & 0.000 & 131.273 & 129.870 & 1.403 & 132.371 & 129.870 & 2.500 \\ \hline
		1.03                & 116.054 & 116.054 & 0.000 & 131.470 & 130.068 & 1.403 & 132.568 & 130.068 & 2.500 \\ \hline
		1.04                & 116.251 & 116.251 & 0.000 & 131.667 & 130.265 & 1.403 & 132.765 & 130.265 & 2.500 \\ \hline
		1.05                & 116.448 & 116.448 & 0.000 & 131.865 & 130.462 & 1.403 & 132.962 & 130.462 & 2.500 \\ \hline
		1.06                & 116.645 & 116.645 & 0.000 & 132.062 & 130.659 & 1.403 & 133.159 & 130.659 & 2.500 \\ \hline
		1.07                & 116.842 & 116.842 & 0.000 & 132.259 & 130.856 & 1.403 & 133.356 & 130.856 & 2.500 \\ \hline
		1.08                & 117.038 & 117.038 & 0.000 & 132.456 & 131.053 & 1.403 & 133.554 & 131.053 & 2.500 \\ \hline
		1.09                & 117.217 & 117.217 & 0.000 & 132.653 & 131.251 & 1.403 & 133.751 & 131.251 & 2.500 \\ \hline
		1.1                 & 117.392 & 117.392 & 0.000 & 132.851 & 131.448 & 1.403 & 133.948 & 131.448 & 2.500 \\ \hline
		1.2                 & 118.913 & 118.913 & 0.000 & 133.672 & 132.259 & 1.413 & 134.760 & 132.259 & 2.500 \\ \hline
		1.3                 & 119.008 & 119.008 & 0.000 & 133.683 & 132.259 & 1.424 & 134.760 & 132.259 & 2.500 \\ \hline
		1.4                 & 119.008 & 119.008 & 0.000 & 133.687 & 132.259 & 1.428 & 134.760 & 132.259 & 2.500 \\ \hline
		1.5                 & 119.008 & 119.008 & 0.000 & 133.687 & 132.259 & 1.428 & 134.760 & 132.259 & 2.500 \\ \hline
		1.6                 & 119.008 & 119.008 & 0.000 & 133.687 & 132.259 & 1.428 & 134.760 & 132.259 & 2.500 \\ \hline
		1.7                 & 119.008 & 119.008 & 0.000 & 133.687 & 132.259 & 1.428 & 134.760 & 132.259 & 2.500 \\ \hline
		1.8                 & 119.008 & 119.008 & 0.000 & 133.687 & 132.259 & 1.428 & 134.760 & 132.259 & 2.500 \\ \hline
		1.9                 & 119.008 & 119.008 & 0.000 & 133.687 & 132.259 & 1.428 & 134.760 & 132.259 & 2.500 \\ \hline
		2                   & 119.008 & 119.008 & 0.000 & 133.687 & 132.259 & 1.428 & 134.760 & 132.259 & 2.500 \\ \hline
		2.5                 & 119.008 & 119.008 & 0.000 & 133.687 & 132.259 & 1.428 & 134.760 & 132.259 & 2.500 \\ \hline
		3                   & 119.008 & 119.008 & 0.000 & 133.687 & 132.259 & 1.428 & 134.760 & 132.259 & 2.500 \\ \hline
		4                   & 119.008 & 119.008 & 0.000 & 133.687 & 132.259 & 1.428 & 134.760 & 132.259 & 2.500 \\ \hline
	\end{tabular}

	\end{adjustbox}
\end{table}

\begin{figure}[tph]
	\begin{center}
		\input{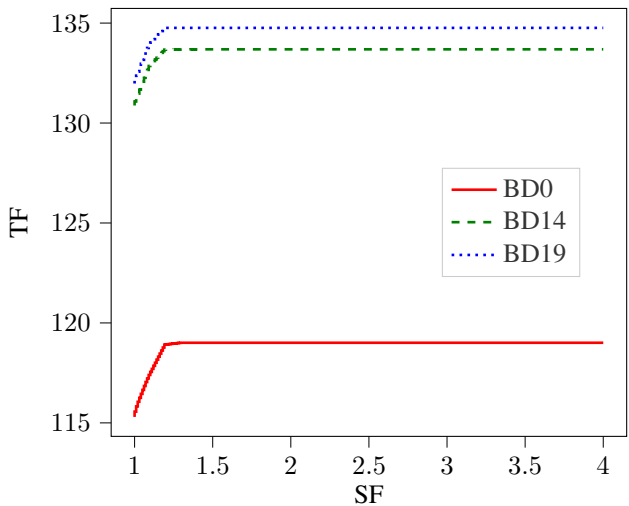}
	\end{center}
	\caption{Trends of TF under Different Transmission Line Capabilities }\label{FIG:lineTF}
\end{figure}

\section{{Conclusions}}\label{sec:conclu}
In this paper, we proposed a real-time flexibility management framework to model economic dispatch with dynamic AGC constraints. The proposed framework can be solved by conducting reformulation and decomposition. We further proposed \lred{seven system flexibility indices} (i.e., TF, EDF, AGCF, EDUPF, EDDNF, AGCUPF, AGCDNF) to reflect the system flexibility and studied how the system factors such that the operational cost budget, ramping capacity and transmission line capacity can impact the system flexibility. We found that the budget of the operational cost can \zb{significantly} contribute to all indicators expect EDDNF, and the improvement of ramping capability can significantly enhance the EDF and the AGCF. In addition, we discover that the transmission capacity only \lgreen{contributes} to EDF.  As the future work, we will study the \zb{power system} flexibility with considering the inertia of the resources in the system.

\bibliographystyle{IEEEtran}
\bibliography{CC_Edge}

\end{document}